\documentclass[12pt]{article}
\usepackage{amsfonts}
\usepackage{amsmath}
\usepackage{amssymb}
\usepackage{color}
\font\smallit=cmti10  

\usepackage{amssymb,latexsym}
\usepackage{latexsym, amssymb, amsmath, amscd, amsfonts, epsfig,graphicx}
\usepackage{pifont}
\usepackage{CJK, psfrag,eepic}
\usepackage{mathrsfs}
\usepackage{txfonts}
\usepackage[perpage,symbol*]{footmisc}
\makeatletter

\renewcommand\section{\@startsection {section}{1}{\z@}
{-30pt \@plus -1ex \@minus -.2ex} {2.3ex \@plus.2ex}
{\normalfont\normalsize\bfseries}}

\renewcommand\subsection{\@startsection{subsection}{2}{\z@}
{-3.25ex\@plus -1ex \@minus -.2ex} {1.5ex \@plus .2ex}
{\normalfont\normalsize\bfseries}}

\renewcommand{\@seccntformat}[1]{\csname the#1\endcsname. }

\makeatother

\newtheorem{thm}{Theorem}[section]

\newtheorem{rem}[thm]{\it Remark}

\newtheorem{lemma}[thm]{Lemma}

 \DeclareMathOperator{\ord}{ord}

 \DeclareMathOperator{\supp}{supp}

\def\qed{\nopagebreak \hfill $\Box$\medbreak}

\begin{document}
\title{
On $n$-sum of an abelian group of order $n$}

%\date{}
\maketitle

\begin{center}
 \vskip 20pt
{\bf Xingwu Xia$^1$ \ Weidong Gao$^2$ }\\ {\smallit $^1$Department
of Mathematics, Luoyang Normal University, Luoyang 471022,P.R. China
}
\\{\smallit $^2$ Center for
Combinatorics, LPMC-TJKLC, Nankai University, Tianjin 300071, P.R.
China}\\

\footnotetext{E-mail address: xxwsjtu@yahoo.com.cn (X.W. Xia),
wdgao1963@yahoo.com.cn (W.D. Gao) }
\end{center}

 \begin{abstract} Let $G$ be an additive finite abelian group of order $n$, and let $S$ be
 a sequence of $n+k$ elements in $G$, where $k\geq 1$. Assume that
 $S$ contains
 $t$ distinct elements. Let $\sum_n(S)$ denote the set consists of all
 elements in $G$ which can be expressed as a sum over
 subsequence of $S$ of length $n$. In this paper we prove that,
 either $0\in \sum_n(S)$ or  $|\sum_n(S)|\geq k+t-1$. This confirms a conjecture
 by Y.O. Hamidoune in 2000.
 \end{abstract}

\section{Introduction}

Let $G$ be an additive  abelian group of order $n$, and let $S=(a_1,
\cdots, a_k)$ be a sequence of elements in $G$ with $k=|S|\geq n$.
Denote by $\sum_n(S)$ the set that consists of all elements which
can be expressed as a sum over a subsequence of $S$ of length $n$,
i.e.
$$
\sum_n(S)=\{a_{i_1}+\cdots+a_{i_n}|1\leq i_1<\cdots
  <i_n\leq k\}.
$$
The famous Erd\H{o}s-Ginzburg-Ziv theorem asserts that if $|S|\geq
2n-1$ then $0\in \sum_n(S)$. The Erd\H{o}s-Ginzburg-Ziv theorem has
attracted a lot of attention and $\sum_n(S)$ has been studied by
many authors (For e.g., see \cite{BL1, BGL, BL, Gao1, Gao2,GL, G1,
Ham, Olson,SC}). In this paper we settle a conjecture by Hamidoune
\cite{Ham} on $\sum_n(S)$
 by showing

%Let $G$ be an additive  abelian group of order $n$, and $S=(a_1,
%\cdots, a_k)$ a sequence of elements in $G$.
%   By $\sigma(S)$ we denote the sum $\sum_{i=1}^ka_i$. For every $1\leq r
%   \leq k$, we define $\sum_r(S)=\{a_{i_1}+\cdots+a_{i_r}|1\leq i_1<\cdots
%   <i_r\leq k\}$, and for every $1\leq l \leq k$, by $\sum_{\geq l}(S)$ we
%   denote the set $\cup_{r=l}^k\sum_r(S)$. We define $h(S)$ to be the
%   maximal number of repeatation of an element in $S$. By supp(S) we denote
%   the set consists of all distinct elements in $S$. If $T$ is a subsequence
%   of $S$, by $ST^{-1}$ we denote the sequence obtained by deleting the
%   terms of $T$ from $S$. Let $A, B$ be two  nonempty subsets
% of $G$, define $A+B=\{a+b|a\in A, b\in B\}$, and  let
% $\mbox{St}(A)=\{h\in G|h+A=A\}$ be the stabilizer of $A$
%   in $G$. Clearly, $\mbox{St}(A)$ is the maximal subgroup $H$ of $G$
%   satisfying $A+H=A$.
%   Y. Hamidoune \cite{Ham} formulated the following conjectures.
%
%
%\medskip
%   {\bf Conjecture A}.
%   {\sl Let $m,n$ be positive integers with $m\geq n+1$, and  Let
%   $S=(x_1, x_2, \cdots, x_m)$ be a sequence of $m$ elements of $Z_n$. Set
%   $k=|\mbox{Supp}(S)|$.
%   Then, one of
%   the following conditions holds:

  % (1) $h(S)\geq n-k+3$.
%
%
%   (2) $\sum_n(S)$ contains a non-null subgroup.
%
%
%   (3) $|\sum_n(S)|\geq m-n+k-1$.}
%
%
%
%\medskip
\begin{thm} \label{theorem1}
Let $G$ be a finite abelian group of order $n$. Let $k\geq 1$ be an
integer. Let
   $S$ be a sequence of $n+k$ elements of $G$. Set
   $t=|\mbox{Supp}(S)|$.
   Then one of the following conditions holds:

   (1) $0\in \sum_n(S)$.

   (2) $|\sum_n(S)|\geq k+t-1$.
\end{thm}

\medskip

Let $G=C_n$ be the cyclic group of order $n$ in Theorem
\ref{theorem1}, we get a positive answer to an open problem by
Hamidoune \cite[Conjecture B]{Ham}.

%We shall obtain Theorem 1 as an easy consequence of the following
%
%
%\begin{th}
%Let $m,n$ be positive integers with $m\geq n+2$, and  Let
%   $S=(x_1, x_2, \cdots, x_m)$ be a sequence of $m$ elements of $Z_n$. Set
%   $k=|\mbox{Supp}(S)|$.
%   Suppose $k\geq 36$. Then, one of
%   the following conditions holds:
%
%   (1) $h(S)\geq n-k+3$.
%
%
%   (2) $\sum_n(S)$ contains a non-null subgroup.
%
%
%   (3) $|\sum_n(S)|\geq m-n+k-1$.
%
%   (4) $0\not\in \sum_n(S)$.
%\end{th}
%
%
%
%
%
%
%From Theorem 1, it is easy to see that, Conjecture B, if true, would
%imply Conjecture A for $k\geq 36$.  If $m\geq n+2$ then Theorem 2 clearly
%implies Theorem 1. For $m=n+1$, clearly, $|\sum_n(S)|=k=m-n+k-1$. So,
%Theorem 2 implies Theorem 1.

\section{Notations and preliminaries}
Let $\mathbb{N}$ denote the set of positive integers, and
$\mathbb{N}_{0}=\mathbb{N}\cup\{0\}$.
 For any two integers $a, b\in \mathbb N_0$, we set
 $[a, b]=\{x \in \mathbb N_0 : a\leq x\leq b\}$.
Throughout this paper, all abelian groups will be written
additively.
 %and for $n, r \in \mathbb N$, we denote by $C_n$ the
%cyclic group of order $n$, and denote by $C_n^r$  the direct sum of
%$r$ copies of $C_n$.

Let $\mathcal F(G)$ be the free abelian monoid, multiplicatively
written, with basis $G$. The elements of $\mathcal F(G)$ are called
sequences over $G$. We write sequences $S\in\mathcal F(G)$ in the
form
$$S=\mathop\Pi\limits_{g\in G}g^{ v_g(S)},\ \text{with}\  v_g(S)\in\mathbb N_0\
\text{for all}\ g\in G.$$
 We call $v_g(G)$ the multiplicity of $g$
in $S$, and we say that $S$ contains $g$ if $v_g(S)>0.$ Further, $S$
is called squarefree if $v_g(S)\leq 1$ for all $g\in G$. The unit
element $1\in\mathcal F(G)$ is called the empty sequence. A sequence
$S_1$ is called a subsequence of $S$ if $S_1\mid S$ in $\mathcal
F(G)$. Let $S_1, \cdots, S_r$ be some subsequences of $S$. We say
$S_1, \cdots, S_r$ are disjoint subsequences if $S_1\cdots S_r|S$.
If a sequence $S\in\mathcal F(G)$ is written in the form
$S=g_1\cdot\ldots\cdot g_l,$ we tacitly assume that $l\in\mathbb
N_0$ and $g_1,\ldots,g_l\in G$.
\par For a sequence
$$S=g_1\cdot\ldots\cdot g_l=\mathop\Pi\limits_{g\in G}g^{v_g(S)}\in\mathcal
F(G),$$ we call
\begin{itemize}
\item $|S|=l=\sum_{g\in G} v_g(G)\in\mathbb N_0$ the $length$ of $S$,
\item $\mathsf h(S)=\max\{\mathsf v_g(S)|g\in G\}\in[0,|S|]\}$ the $maximum$ of
the $multiplicities$ of $S$,
\item$\supp(S)=\{g\in G| v_g(S)>0\}\subset G$ the $support$ of $S$,
\item$\sigma(S)=\sum_{i=1}^lg_i=\sum_{g\in G}v_g(S)g\in G$ the $sum$ of $S$,
\end{itemize}

\noindent The sequence $S$ is called
\begin{itemize}
\item a $zero-sum\  sequence$ if $\sigma(S)=0$,
\item $zero-sumfree$ if $0\not\in\sum(S)$,
\end{itemize}

\noindent For every $1\leq r
   \leq k$, we define $$\sum_r(S)=\{a_{i_1}+\cdots+a_{i_r}|1\leq i_1<\cdots
   <i_r\leq k\}.$$ For every $\ell \in [1,k]$, define  $$\sum_{\geq \ell}(S)=\bigcup_{r=\ell}^k\sum_{\ell}(S)$$ and
   $$\sum_{\leq \ell}(S)=\bigcup_{r=l}^\ell\sum_{\ell}(S).$$ Define
$$\sum(S)=\bigcup_{r=1}^k\sum_r(S).$$
  %Denote by $h(S)$  the
%   maximal multiplicity of $S$, i.e.
%   $$
%   h(S)=\max_{g\in G}\{v_g(S)\}.
%   $$

\noindent Let $A,B$ be two nonempty subsets of $G$. Define
$$
A+B=\{a+b: a\in A, b\in B\}.
$$
If $A=\{x\}$ for some $x\in G$ then we simply denote $A+B$ by $x+B$.
For any nonempty subset $C$ of $G$, let $-C=\{-c:c\in C\}$. For
every $g\in G$, let $\gamma_g(A,B)$ denote the number of the pairs
of $(a,b)$ such that $a\in A,b\in B$ and $a+b=g$.

  \bigskip
\noindent To prove Theorem \ref{theorem1} we need some preliminaries
begin with

%\begin{lemma}(see \cite{Kn} or \cite{Na})
%     Let $G$ be an abelian group, and let $A_1, \cdots, A_r$
%   be some nonempty subsets of $G$. Then,
%   $|A_1+\cdots +A_r|\geq |A_1+H|+\cdots +|A_r+H|-(r-1)|H|$, where
%   $H=\mbox{St}(A_1+\cdots +A_r)$.
%\end{lemma}

   \begin{lemma}\cite{S}\label{ps}
    Let $G$ be a abelian group, and let $A, B$ be two finite
   subsets of $G$ with $ A\cap (-B)=\{0\}$.  Then, $|A+B|\geq |A|+|B|-1$.
\end{lemma}

\noindent   By using Lemma \ref{ps}  repeatedly one can prove

   %\begin{lemma}\cite{Gao1}\label{wdg}
%   Let $G$ be a finite abelian group of order $n$, and let $S$ be a
%   sequence of $n$ elements in $G$. Set $h=h(S)$. Then,
%   $0\in \sum_{\leq h}(S)$.
%\end{lemma}

\begin{lemma} \label{ben} \cite{BEN} Let $S$ be a zero-sum free
sequence over an abelian group, let $S_1, \cdots, S_k$ be some
disjoint subsequences of $S$. Then,
$$
\left|\sum(S)\right|\geq \sum_{i=1}^k \left|\sum(S_i)\right|.
$$
\end{lemma}

\noindent   The following lemma is crucial in this paper.

   \begin{lemma} \cite{Gao1} \label{gwd}
   Let $G$ be a finite abelian group of order $n$, and let
   $S=0^mT \in \mathcal F(G)$ be a sequence of
   length $|S|\geq n$.   If $h(T)\leq m$ then
   $$\sum_{\geq n-m}(T)=   \sum_n (S).$$
\end{lemma}

\begin{rem} Lemma \ref{gwd} in the present version appeared first in \cite{BGL} and then in \cite{SC}.
The proof of \cite[Theorem 3]{Gao1} implies that Lemma \ref{gwd}.
\end{rem}

%\proof  The method used in \cite{Gao1} works for this theorem. For the proof
%is not long and for the completeness we include the  proof here.
%
%Clearly, $\sum_n (S)\subset \sum_{\geq n-k}(S)$.
%   So it suffices to prove that $\sum_{\geq n-k}(S)\subset \sum_n (S)$. Let
%   $x\in \sum_{\geq n-k}(S)$, we have to show $x\in \sum_n (S)$. By the
%   definition of $\sum_{\geq n-k}(S)$, $x=\sigma (W)$ for some subsequence
%   $W$ of $T$ with $|W|\geq n-k$. If $|W|\leq
%   n$, then $W'=W0^{n-|W|}$ is a subsequence of $S$ with
%   $|W'|=n$ and $x=\sigma(W')\in \sum_n(S)$. Otherwise,
%   $|W|\geq n+1$. Since $h=h(W)\leq h(T)\leq k$, by using
%   Lemma \ref{wdg} repeatedly, one can find some disjoint zero-sum subsequences
%   $W_1, \cdots, W_l$ of $W$ such that  $|W(W_1\cdots W_l)^{-1}|\leq n-1$ ,
%   $|W(W_1\cdots W_{l-1})^{-1}|\geq n$ and $|W_i|\leq k$ for $i=1, \cdots,
%   l$. Therefore, $n-1\geq |W(W_1\cdots W_l)^{-1}|= |W(W_1\cdots
%   W_{l-1})^{-1}|-|W_l|\geq n-k$.
%   Notice that $\sigma(W)=\sigma(W(W_1\cdots W_l)^{-1})$ and it reduces to
%   the case above. \qfd

  \begin{lemma} \cite{EE}\label{ee}
   Let $S$ be a  subset of an abelian group $G$
   with $0\not\in \sum(S)$. Then,
\begin{enumerate}
\item $|\sum(S)|\geq 2|S|-1$;

\item if $|S|\geq 4$ then $|\sum(S)|\geq
   2|S|$;

\item  if $|S|= 3$ and $S$ does not contain exactly one element of order two then $|\sum(S)|\geq
   2|S|$.
\end{enumerate}
\end{lemma}

{\it Proof}. 1. and 2. has been proved in \cite{EE}.

\medskip
3. If $S$ contains no element of order two, then the result has been
proved also in \cite{EE}. Now assume that $S$ contains at least two
elements of order two. Let $S=\{a,b,c\}$ with $\ord(a)=\ord(b)=2$.
If $a+b=c$ then $a+b+c=0$, a contradiction. Therefore, $a+b\not \in
S$. If $a+c=b$ then $a+c+b=0$, also a contradiction.  Hence,
$a+c\not\in S$. Similarly, we can prove $b+c\not \in S$. Note that
$a+b+c\not \in \{a,b,c,a+b,b+c,c+a\}$. Therefore, $|\sum(S)|=7$ and
we are done. \qed

\begin{lemma}  \label{a.g}
 Let $G$ be a finite abelian group, and let
   $S\in \mathcal F(G)$  be a zero-sumfree sequence.
   Then, $\sum(S)\geq|S|+|\mbox{Supp}(S)|-1$, and we have strict inequality
   except $|S|\leq 2$, or $|S|=3$ and $S$ contains exactly one element of
   order two.
\end{lemma}

{\it Proof}.  Let $S_1$ be a squarefree subsequence of $S$ with
$|S_1|=|\supp(S)|$, and let $S_2=SS_1^{-1}$. Apply Lemma \ref{ben}
to $S_2$  we obtain that
$$
|\sum (S_2)|\geq |S_2|.
$$
Again apply Lemma \ref{ben} to $S=S_1S_2$ we obtain that
$$
|\sum(S)|\geq |\sum(S_1)|+|\sum(S_2)|\geq |S_2|+|\sum(S_1)|=
|S|-|S_1|+|\sum(S_1)|.
$$
Now the result follows from Lemma \ref{ee}. \qed

\begin{lemma} \label{lemma8.9} Let $A$ be a finite subset of an abelian
group  with $A\cap -A=\{0\}$. If $|A|\in [3,5]$ then
$|A\dot{+}A|\geq |A|.$
\end{lemma}

{\it Proof}. Assume to the contrary that $|A\dot{+}A|\leq |A|-1.$
Since $0\in A$ we infer that $A\setminus \{0\}\subset A\dot{+}A$. It
follows that $A\setminus \{0\}=A\dot{+}A.$

 Let $x\in A\setminus \{0\}$. Then, $A\setminus \{0\}=A\dot{+}A$
 implies that
\begin{equation} \label{eq8.0}
A\setminus \{0,x\}=x+(A\setminus \{0,x\})
\end{equation}
holds for every $x\in A \setminus \{0\}.$

 \noindent Therefore,
 $\sum_{y\in A\setminus \{0,x\}}y=\sum_{y\in A\setminus
 \{0,x\}}(x+y)$. It follows that
 \begin{equation} \label{eq8.1}
 (|A|-2)x=0
 \end{equation}
 holds for every $x\in A.$

\noindent Equality (\ref{eq8.1}) implies that $|A|-2\geq 2$. Hence,
 $$
 |A|\in [4,5].
 $$
 If $|A|=4$ then (\ref{eq8.1}) gives that $2x=0$ for every $x\in
 A\setminus \{0\}.$ Thus, $A\cap (-A)=A$, a contradiction. So, we may
 assume that $$|A|=5.$$
 Let $A=\{0,a,b,c,d\}$. Now (\ref{eq8.1}) gives that
 $$
 3a=3b=3c=3d=0.
 $$
From (\ref{eq8.0}) we may assume that
$$
a+b=c
$$
and
$$
\{a+c,a+d\}= \{b,d\}.
$$
If $a+c=b$ then $a+b=c$ gives that $2a=0$. This together with $3a=0$
implies that $a=0$, a contradiction. Therefore,
$$
a+c=d \mbox{ and } a+d=b.
$$
Now we have $b+c=a+d+a+b=b+d-a$. Hence,
$$
b+c=b+d-a.
$$
This implies that $b+c\not \in \{b,c,d\}$. It follows from
(\ref{eq8.0}) that $b+c=b+d-a=a$. Therefore, $2a=b+d$. Again by
(\ref{eq8.0}) we obtain that $-a=2a=b+d\in A \setminus \{0\}\subset
A.$ Hence, $a\in A\cap (-A),$ a contradiction.
 \qed

\begin{lemma} \label{lemma9.0} Let $A$ be a finite subset of an abelian
group $G$. Suppose that $A\cap -A=\{0\}$ and suppose that $|A|\geq
6$. Then, $|A\dot{+}A|\geq |A|+1.$
\end{lemma}

{\it Proof.} Assume to the contrary that $|A\dot{+}A|\leq |A|.$ For
every $x\in A\setminus \{0\}$, let $$A_x=A\setminus \{0,x\}.$$ Since
$0\in A$, we infer that $A\setminus \{0\} \subset A\dot{+}A$. It
follows from $|A\dot{+}A|\leq |A|$ that $|A\dot{+}A\setminus A|\leq
1$. Therefore, $|(x+A_x)\setminus A|\leq 1$. By the hypothesis that
$A\cap -A=\{0\}$ we deduce that $(x+A_x)\setminus A=(x+A_x)\setminus
A_x$. Hence, $|(x+A_x)\setminus A_x|\leq 1$, and this is equivalent
to

\begin{equation} \label{eq00}
|(x+A_x)\cap A_x|\geq |A_x|-1.
\end{equation}
 Note that $|(x+A_x)\cap
A_x|=|(-x+A_x)\cap A_x|$. Thus,

\begin{equation} \label{eq01}
|(-x+A_x)\cap A_x|\geq |A_x|-1.
\end{equation}

 We
assert that one of the following statements hold:
\begin{enumerate}
\item $2x=y+z$ for some $y,z \in A_x$ with $y\neq z$.

\item $2x=2y$ for some $y\in A_x$.

\item $x=2y$ for some $y\in A_x$.

\item $x\not \in A_x+A_x$.
\end{enumerate}

If $|(x-A_x)\cap A_x|\geq 2$  then  by (\ref{eq01}) we have
$(-x+A_x)\cap (x-A_x) \neq \emptyset$, and therefore 1. or 2. holds.

Otherwise, $|(x-A_x)\cap A_x|\leq 1$. If $|(x-A_x)\cap A_x|= 1$ then
we must have 3. holds. For the remainder case that $|(x-A_x)\cap
A_x|=0$, we get 4. holds. This proves the assertion.

For every $i\in [1,4]$, let $B_i$ be the subset of $A$ consisting of
all elements $x\in A\setminus \{0\}$ such that the item i. in the
assertion holds for $x$. Then,
\begin{equation} \label{eq02}
|B_1|+|B_2|+|B_3|+|B_4|\geq |A|-1.
\end{equation}

Let $C=\{2x: x\in A\} \setminus (A\dot{+}A\}$. Then, $0\in C$ and
$A\dot{+}A=(A+A)\setminus C$. By Lemma \ref{ps}, $|A+A|\geq 2|A|-1$.
It follows from $|A\dot{+}A|\leq |A|$ that

\begin{equation} \label{eq03}
|C|\geq |A|-1.
\end{equation}
From (\ref{eq03}) we deduce that $$|B_1|\leq 1.$$
 Note that
$A\setminus \{0\} \subset A\dot{+}A$ and again from (\ref{eq03}) we
deduce that.
$$
|B_3|\leq 1.
$$
We show next that
$$
|B_4|\leq 1.
$$
Assume to the contrary that $|B_4|\geq 2.$ Let $x,y\in B_4$ with
$x\neq y$, and let $z\in A\setminus \{0,x,y\}$. Then, $x,y,z \not\in
z+A_z\cap A_z$. Hence, $|z+A_z\cap A_z|\leq |A_z|-2$, a
contradiction on (\ref{eq00}). Now by (\ref{eq02}) we infer that
$|B_1|+|B_2|\geq |A|-1-2\geq 3.$ Since $|B_1|\leq 1$, we infer that
$|B_2|\geq 3$, or $|B_1|=1$ and $|B_2|\geq 2$. But in both cases we
have $|C|\leq |A|-2$, a contradiction on (\ref{eq03}). \qed

\section{Proof of Theorem \ref{theorem1}}

\medskip
 {\it  Proof of Theorem \ref{theorem1}}. Without loss of
   generality, we may assume that $S=0^{h}T$ with
   $0\not |T$  and $h=h(S)$. Then $|T|=n-h+k$.
Assume that $0\not \in \sum_n(S)$, we need to show
$$|\sum_{n}(S)|\geq k+t-1.$$ By Lemma \ref{gwd}, it suffices to prove that
$$
|\sum_{\geq n-h}(T)|\geq k+t-1.
$$

    Suppose that ${\rm supp}(T)=\{x_1, \cdots, x_{t-1}\}$.

 Let $T_{0}$ be one of the maximal (in
   length)  subsequence of $T$ with $\sigma(T_{0})=0$ ($T_0$ is the empty sequence if $T$ is zero-sum free).
   By renumbering if necessary we assume that
   $$\mbox{supp}(T_{0})=\{x_1, \cdots, x_r\}$$ for some $r\in
   [1,t-1].$ By Lemma \ref{gwd},
   $$
   |T_0|\leq n-h-1.
   $$

   Let
   $$
   T_1=TT_0^{-1}.
   $$
   Then, $T_1$ is zero-sum free.
   %We distinguish two cases:
%
%   {\bf Case 1.} $\mbox{supp} (T_0) \not \subseteq \mbox{supp}
%   (T_1)$.
By renumbering we may assume  that
$$\mbox{supp} (T_0) \setminus \mbox{supp}
   (T_1)=\{x_1, \cdots, x_{\ell}\}$$ for some $\ell \in [0,r]$. It
   follows that
$$
\mbox{supp}
   (T_1)=\{x_{\ell+1}, \cdots, x_{t-1}\}.
$$

{\bf Claim 1.} $\{x_1, \cdots,x_{\ell}\}\cap \sum(T_1)=\emptyset$.

\medskip
Assume to the contrary that $x_i=\sigma(V_1)$ for some $i\in
[1,\ell]$ and $V_1|T_1$. By the definition of $\{x_1,
\cdots,x_{\ell}\}$ we deduce that $|V_1|\geq 2$. Therefore,
$T_0x_i^{-1}V_1$ is a zero-sum subsequence of $T$ of length
$|T_0|-1+|V_1|>|T_0|$, a contradiction with the maximality of $T_0$.
This proves Claim 1.

\medskip
   Note that $|T_1|=n-h+k-|T_0|\geq k+1$. We can choose a
   subsequence $V$ of $T_1$ with $|V|=k-1$. Let $U=T_1V^{-1}$.
   Then,
   $$
   |U|=n-h-|T_0|+1.
   $$

   We choose $V$ so that

   (i) $|\mbox{supp}(V)|$ attain the maximal value among all subsequences of $T_1$ of length
   $k-1$;

   (ii) $|\mbox{supp}(V)\cap \mbox{supp}(U)|$ attains the maximal value
   subject to (i).

Let
$$
A=\{0,-x_1,\cdots,-x_{\ell}\}.
$$
Clearly,
$$A\subset \sum_{\geq |T_0|-1}(T_0).$$
Let
$$
B=\{\sigma(U)\}\bigcup(\sigma(U)+\sum(V)).
$$
Clearly,
\begin{equation} \label{eq1}
|A|=\ell+1 \mbox{ and } |B|=1+|\sum(V)|
\end{equation}
and
\begin{equation} \label{eq2}
A+B\subset \sum_{\geq n-h}(T).
\end{equation}

By Lemma \ref{a.g} , we have that
\begin{equation} \label{eq3}
 |B|=1+|\sum(V)|\geq |V|+|\mbox{supp}(V)|=k-1+|\mbox{supp}(V)|.
\end{equation}

 Let
$$
C=\{\sigma(U)-x: x\in \mbox{supp}(U) \}.
$$
Then,
\begin{equation} \label{eq4}
|C|=|\mbox{supp}(U)|
\end{equation}
Since $\sigma(U)-x=\sigma(T_0U)-x$, we infer that
\begin{equation} \label{eq5}
C\subset \sum_{\geq n-h}(T).
\end{equation}

\medskip
{\bf \mbox{Claim 2.}} $|A+B|\geq |A|+|B|-1.$

\medskip
If $\gamma_{\sigma(U)}(A,B)>1$, then we deduce that
$\sigma(U)=-x_{i}+(\sigma(U)+\sigma(V_{1}))$ for some $i\in
[1,\ell]$ and some subsequence $V_{1}$ of $V$.  It follows that
$x_{i}=\sigma(V_{1}))$, a contradiction with Claim 1. Therefore,
$$\gamma_{\sigma(U)}(A,B)=1.$$ Let $B'=-\sigma(U)+B$. Then,
$\gamma_{\sigma(U)}(A,B)=1$ implies that $A\cap (-B')=0$. It follows
from Lemma \ref{ps} that
$$
|A+B|=|A+B'|\geq |A|+|B'|-1=|A|+|B|-1.
$$

\medskip

${\bf Claim 3.} (A+B)\cap C=\emptyset.$

Assume to the contrary that Claim 3 is false. It follows that we
have the following possibilities:

(a) $\sigma(U)-x=\sigma(U)$ with $x\in \mbox{supp}(U)$;

(b) $\sigma(U)-x=\sigma(U)+\sigma(V_1)$ with $x\in \mbox{supp}(U)$
and $V_1|V$;

(c) $\sigma(U)-x=\sigma(U)-x_i$ with $x\in \mbox{supp}(U)$ and $i\in
[1,\ell]$;

(d) $\sigma(U)-x=\sigma(U)-x_i+\sigma(V_1)$ with $x\in
\mbox{supp}(U)$, $i\in [1,\ell]$ and $V_1|V$.

Possibility (a) implies that $x=0$, a contradiction; Possibility (b)
implies that $\sigma(xV_1)=0$, a contradiction on $T_1$ is zero-sum
free; Possibility (c) implies that $x=x_i$, a contradiction on the
definition of $A$; and Possibility (d) implies that
$x_i=\sigma(xV_1)$,  a contradiction on Claim 1. This proves Claim
3. Now from Lemma \ref{a.g} , Claim 2, Claim 3, and equation (1-5)
we obtain that
$$
\begin{array}{ll} |\sum_{\geq n-h}(T)| & \geq |A+B|+|C| \\ & =|A+B|+|\mbox{supp}(U)| \\ & \geq
|A|+|B|-1+|\mbox{supp}(U)| \\ & =1+\ell+|\sum(V)|+|\mbox{supp}(U)|
\\ & \geq 1+\ell+|V|+ |\mbox{supp}(V)|-1+ |\mbox{supp}(U)|\\
&=\ell+k-1+|\mbox{supp}(V)|+ |\mbox{supp}(U)|
\\ & =\ell+k-1+|\mbox{supp}(UV)|+|\mbox{supp}(U\cap V)| \\ &
=k+t-2+|\mbox{supp}(U\cap V)|. \end{array}
$$
This forces that $\mbox{supp}(U\cap V)=\emptyset$, $|A+B|=|A|+|B|-1$
and $|\sum(V)|=|V|+|\mbox{supp}(V)|-1$. It follows from the choice
of $U$ and $V$ that

$$T_1=UV \mbox{ is square free}.
$$
Now by Lemma \ref{a.g} and $|\sum(V)|=|V|+|\mbox{supp}(V)|-1$ we
infer that $k-1=|V|\leq 3.$ If $k-1=3$ then $|T_1|\geq k+1=5$. Since
$T_1$,  we can choose $V$ so that $V$ contains no element with order
two, or $V$ contains at least two elements with order two. Now again
by Lemma \ref{a.g} we have that $|\sum(V)|\geq |V|+|\mbox{supp}(V)|$
and repeat the process above we obtain that $|\sum_{\geq n-h}(T)|
\geq |A+B|+|C|\geq k+t-1$, a contradiction. Therefore,
$$
k\leq 3.
$$
Note that $|\sum_{\geq n-h}(T)|=|\sigma(T)-\{0\}\cup\sum_{\leq
k}(T)|=|\{0\}\cup\sum_{\leq k}(T)|.$ It suffices to prove that
\begin{equation} \label{eq6}
|\{0\}\cup\sum_{\leq k}(T)|\geq t+k-1.
\end{equation}
Recall that $T_1$ is square free. Let $D=\{0\}\cup T_1$.

\medskip
 If $k=1$ then $|\{0\}\cup\sum_{\leq 1}(T)|=1+|\mbox{supp}
(T)|=t= t+k-1$ and we are done.

\smallskip
If $k=2$ then $|T_1|\geq k+1\geq 3.$ By Claim 1, Lemma
\ref{lemma8.9} and Lemma \ref{9.0} we obtain that
$|\{0\}\cup\sum_{\leq 2}(T)|\geq 1+\ell+|\sum_{\leq
2}(T_1)|=1+\ell+|D\dot{+}D|\geq 1+\ell+|T_1|+1=t+1$ and we are done.

\smallskip
So it remains to consider the case that $$k=3.$$ Now we have
$$
|D|=1+|T_1|\geq 2+k\geq 5.
$$
If $|D|=5$ then $|T_1|=4$. By Lemma \ref{ee} we have that
$|\sum_{\leq 3}(T_1)|= |\sum(T_1)|-1 \geq 2|T_1|-1.$ It follows from
Claim 1 that $|\{0\}\cup\sum_{\leq 3}(T)|\geq 1+\ell+|\sum_{\leq
3}(T_1)|\geq 1+\ell+2|T_1|-1=t+|T_1|-1>t+k-1$. This proves
(\ref{eq6}) for this case.

\smallskip
Now assume that $|D|\geq 6.$ By Lemma \ref{lemma9.0} we have that
$|\sum_{\leq 2}(T_1)|=|D\dot{+}D|\geq |D|+1=|T_1|+2.$ It follows
from Claim 1 that $|\{0\}\cup\sum_{\leq k}(T)| \geq
|\{0\}\cup\sum_{\leq 2}(T)|\geq 1+\ell+|\sum_{\leq 2}(T_1)|\geq
1+\ell+|T_1|+2=t+2=t+k-1.$ This proves (\ref{eq6}) and completes the
proof. \qed

\medskip
{\bf Acknowledgements.}  This work was supported  by the National
Key Basic Research Program of China (Grant No.  2013CB834204), the
PCSIRT Project of the Ministry of Science and Technology,  the
National Science Foundation of China, and the Education Department
of Henan Province (Grant No. 2009A110012)


\begin{thebibliography}{99}
\bibitem{BL1} B. Bollob\'{a}s and I. Leader, \emph{The number of $k$ -sums modulo $k$},
J. Number Theory 78(1999)27-35.

\bibitem{BGL} A. Bialostocki, D. Grynkiewicz and M. Lotspeich  \emph{On some developments of the Erd?s-Ginzburg-Ziv
theorem. II}, Acta Arith. 110(2003)173-184.

\bibitem{BL} A. Bialostocki and M. Lotspeich, {\it Some developments of
the Erd\"{o}s-Ginzburg-Ziv theorem}, in "Sets , Graphs and
Numbers", Coll. Math. Soc. J.Bolyai, 60(1992), 97-117.

\bibitem{BEN} J.D. Bovey, P. Erd\H{o}s, I. Niven, \emph{Conditions for zero sum modulo n},
Canad. Math. Bull. 18(1975)27-29.

\bibitem{Caro1} Y. Caro, {\it Zero-sum problems- a survey}, Discrete
Math., 152(1996), 93-113.

\bibitem{EE} R.B. Eggleton and P. Erd\H{o}s, {\it two combinatorial
problems in group theory}, Acta Arith. 21(1972) 111-116.


\bibitem{Gao1} W.D. Gao, {\it A combinatorial problem on finite abelian
groups,} J. Number Theory 58(1996) 100-103.

\bibitem{Gao2} W.D. Gao, \emph{On the number of subsequences with given
sum}, Discrete Mathematics, 195(1999) 127-138.

\bibitem{GL} W.D. Gao and I. Leader, \emph{ Sums and $k$ -sums in abelian groups of order $k$},
 J. Number Theory 120(2006) 26-32.


\bibitem{GH1}A.Geroldinger and F.Halter-Koch, {\it Non-Unique Factorizations.
Algebra, Combinatorial and Analytic Theory, Pure and Applied
Mathematics},vol.278, Chapman $\&$ Hall/CRC,2006.

\bibitem{G1} D. Grynkiewicz, \emph{On a conjecture of Hamidoune for subsequence
sums}, Integers, 5(2005) No.2, A7.

\bibitem{Ham} Y. O. Hamidoune, {\it Subsequence sums},
   Combinatorics, Probability and Computing 12(2003) 413-425

%\bibitem{Kn} M. Kneser, {\it Ein satz \"{u}ber abelsche gruppen mit
%anwendungen auf die geometrie der zahlen}, Math. Z., 64(1955), 429-434.

\bibitem{S}  P. Scherk, {\it Distinct elements in a set of sums},
Amer. Math. Monthly, 62(1955), 46-47.

 \bibitem{Na} M. B. Nathanson, {\it Additive Number Theory, Inverse Problems
 and the Geometry of Sumsets}, GMT 165, Springer-Verlag (New York, 1996).

 \bibitem{Olson} J.E. Olson, \emph{An addition theorem for finite abelian
 groups}, J. Number Theory 9(1977) 63-70.



\bibitem{SC} S. Savchev and F. Chen, \emph{Long n -zero-free sequences in finite cyclic groups},
 Discrete Mathematics, 308(2008)1¨C8.
\end{thebibliography}
\end{document}